\documentclass[12pt]{article}

\usepackage{graphicx,amssymb}
\usepackage{graphicx}
\usepackage{dcolumn}
\usepackage{amsfonts}
\usepackage{amssymb}
\usepackage{amsmath}
\usepackage{latexsym}
\usepackage[mathscr]{eucal}
\usepackage{epsfig}

\newif{\ifcomentarios}
\comentariosfalse

\newtheorem{theorem}{Theorem}

\newtheorem{definition}[theorem]{Definition}

\newcommand{\be}{\begin{eqnarray}}
\newcommand{\en}{\end{eqnarray}}
\newcommand{\bee}{\begin{eqnarray*}}
\newcommand{\ene}{\end{eqnarray*}}

\newtheorem{Teorema}{Theorem}
\newtheorem{Lema}{Lemma}
\newtheorem{Propo}{Proposition}

\newtheorem{Obs}{Remark}
\renewcommand{\mathbf}{\boldsymbol}

\begin{document}

\title{Conformal Deformation from Normal to Hermitian Random Matrix Ensembles%
}
\author{A.M. Veneziani$^1$, T. Pereira$^2$, and D.H.U. Marchetti$^1$ \\
$^1$ Universidade de S\~ao Paulo, Instituto de F\'{\i}sica \\
$^2$ Universidade Federal do ABC, Centro de Matem\'atica\\
E-mail: alexei@if.usp.br; tiago.pereira@ufabc.edu.br;
marchetti@if.usp.br}
\date{}

\maketitle

\begin{abstract}
We investigate the eigenvalues statistics of ensembles of normal random
matrices when their order $N$ tends to infinite. In the model the
eigenvalues have uniform density within a region determined by a simple
analytic polynomial curve. We study the conformal deformations of normal
random ensembles to Hermitian random ensembles and give sufficient
conditions for the latter to be a Wigner ensemble.
\end{abstract}

\section{Introduction and Statement of Results}

Since early fifties Hermitian random matrix theory plays an important role
in the statistical description of the spectra of complex systems \cite%
{Mehta,Deift}. Recently non-Hermitian random matrices have been used to
treat problems in superconductor physics with columnar defects \cite%
{Efetov,Hatano}, in quantum chaotic systems \cite{ChaosR}, and in quantum
chromodynamics \cite{QCD,Sommers2003}.

Normal random matrix ensembles have been playing a major role in several
areas such as in the study of fractional quantum Hall effect \cite{ChauPLA},
quantum Hele-Shaw flows \cite{Hedenmalm}, integrable hierarchies \cite%
{Wiegmann}, and integrable structure of the Dirichlet boundary problem \cite%
{EFensemble1,EFensemble2}.

In the present work the ability of normal ensembles to be conformally
deformed into Hermitian ensembles is exploited to project density of
eigenvalues of non-Hermitian matrices into the real axis. The problem is
addressed by using the so-called \textit{invariant ensemble model},
characterized by the probability of finding a $N\times N$ matrix $M$ of a
class within the ensemble given by 
\begin{equation}
P(M)dM\propto \exp \left\{ -N\text{Tr}[V(M)]\right\} dM,  \label{Pm}
\end{equation}%
with the trace $\text{Tr}[V(M)]$ and the Riemann volume $dM$ invariant under
unitary transformations. The corresponding eigenvalue density, in the limit $%
N\rightarrow \infty $, depends on the particular form of $V(M)$. For the
Wigner ensemble of Hermitian matrices with $V(M)=\dfrac{1}{\sigma ^{2}}%
M^{\ast }M$ ($M^{\ast }$ is the Hermitian conjugate of $M$), the entrances
of $M$ are independent and identically distributed Gaussian random variables
with zero mean and variance $\sigma ^{2}/N$. The density of eigenvalues
follows the Wigner semicircle law supported on $[-2\sigma ,2\sigma ]$ \cite%
{Deift}: 
\begin{equation}
d\mu _{W}(x)=\frac{1}{2\pi \sigma }\sqrt{4-x^{2}/\sigma ^{2}}\chi _{\lbrack
-2\sigma ,2\sigma ]}\left( x\right) dx,  \label{miW}
\end{equation}%
where $\chi_{A}(x)=1$ if $x \in A$ and $0$ otherwise.

A particularly interesting potential has been put forward by Wiegmann,
Zabrodin and coworkers \cite{Wiegmann,EFensemble1,EFensemble2,EFensemble3}
who established a connection between normal random matrices and conformal
mappings. They considered 
\begin{equation}
V(M)=\frac{1}{t_{0}}(M^{\ast }M-p(M)-p(M)^{\ast }),  \label{FEF}
\end{equation}%
where%
\begin{equation}
p(z)=\sum_{j\geq 1}t_{j}z^{j}  \label{pz}
\end{equation}%
with $t_{0}>0$ and $t_{j}\in \mathbb{C}$. As $N \rightarrow \infty$, they
showed, at the level of formal manipulations, that : $(A)$ the density of
eigenvalues is uniform within a simply connected domain $D\subset \mathbb{C}$
whose boundary is given by a simple analytic curve $\gamma $; $(B)$ the
domain $D$ is characterized by the fact that its exterior harmonic moments 
\begin{equation}
t_{j}=\frac{1}{2\pi ij}\oint_{\gamma }\bar{z}z^{-j}dz,\qquad j\geq 1~,
\label{tj}
\end{equation}%
where $\pi t_{0}$ stands for the area of $D$, are the coefficients of (\ref%
{pz}); and $(C)$ the Riemann mapping from the exterior of the unit disk onto
the exterior of the domain $D$ obeys, as a function of the $t_{j}$, the
equations of the integrable dispersionless Toda hierarchy.

Potentials of the form (\ref{FEF}) give rise two sorts of mathematical
problems. Except in the case of polynomial $p(z)$ of degree $2$, where the
domain $D$ is bounded by an ellipse, $V(z)$ is not bounded from bellow and
integrals with respect to (\ref{Pm}) diverge. The other problem concerns
with the fact that $D$ may not be uniquely determined by the moments (\ref%
{tj}). From the point of view of equilibrium measures (see Section \ref{Pre}%
), a relevant fraction of eigenvalues of a $M$ may escape to infinity or to
another Riemann surface.

Recently, the results $(A)$ and $(B)$ have been set in a rigorous frame by
Elbau and Felder \cite{ElbauFelder}. To avoid the above mentioned problems,
they consider the following restrictions: \medskip \noindent 

\noindent \textbf{Elbau-Felder Potential.} If%
\begin{equation}
p(z)=t_{1}z+t_{2}z^{2}+\cdots +t_{n+1}z^{n+1}~  \label{pzn}
\end{equation}%
is an analytic polynomial of degree $n+1$ with $t_{0}>0$ and $\mathbf{t}%
=\left( t_{1},\ldots ,t_{n+1}\right) \in \mathbb{C}^{n+1}$ such that $t_{1}=0
$, $\left\vert t_{2}\right\vert <1/2$, Elbau--Felder potential is a
real--valued function on $\mathbb{C}$ given by%
\begin{equation*}
V(z)=\frac{1}{t_{0}}(\left\vert z\right\vert ^{2}-p(z)-\overline{p(z)})~.
\end{equation*}

It can be shown by direct computations that $V(z)$, under the above
conditions, is positive in a neighborhood of $z=0$ and has a non-degenerate
absolute minimum at $z=0$. From now on, $V(z)$ shall stand for the
Elbau-Felder potential. The problem of divergence of the integrals is solved
by Elbau--Felder in a na\"{\i}ve way -- imposing that the eigenvalues of
matrices within the normal ensemble remains bounded: \medskip

\noindent \textbf{Elbau-Felder Ensemble.} Let $\Sigma \subset \mathbb{C}$ be
the closure of a bounded open set that contains the origin and consider the
following class of matrices 
\begin{equation}
\mathcal{N}_{N}(\Sigma )=\{A\in \text{Mat}_{\mathbb{C}}(N):[A,A^{\ast
}]=0,\sigma (A)\subset \Sigma \}  \label{c}
\end{equation}%
where $\sigma (A)$ denotes the spectrum of $A$. An ensemble is said to be of
Elbau-Felder type of degree $n+1$ if it fulfills conditions stated between (%
\ref{pzn}) and (\ref{c}). A closed polynomial curve $\gamma $ of degree $n$
can be parametrized by%
\begin{equation}
w\mapsto h(w)=rw+\sum_{j=0}^{n}a_{j}w^{-j}\ ,\qquad \left\vert w\right\vert
=1  \label{h}
\end{equation}%
for some $r>0$ and the $a_{j}\in \mathbb{C}$. Elbau and Felder have shown
that, as long as $\left\vert t_{2}\right\vert <1/2$ and $t_{0}$ is small
enough, the problem of determining the exterior moments $t_{j}$ out of the
curve has a unique solution for simple closed analytic polynomial curves.
They give a set of equations that defines an invertible map $%
F:(r^{2},a_{0},\ldots ,a_{n}/r^{n})\longrightarrow (t_{0},\ldots ,t_{n+1})$
from $\mathbb{R}\times \mathbb{C}^{n+1}$ into itself about $(0,0,2\bar{t}%
_{2},\ldots ,(n+1)\bar{t}_{n+1})$ ($\bar{t}_{j}$ stand for the complex
conjugate of the $t_{j}$). By the Euler--Lagrange variational equations, the
eigenvalues density is uniform in $D$. We refer to Theorem \ref{ElbauFelder}
for a precise statement.

In the present work, we study conformal deformations of the Elbau-Felder
ensembles into Hermitian ensembles. This is achieved by considering a family
of polynomial curves of degree $n$: $w\mapsto h(w;s)$, with the $a_{j}(s)$
depending on a parameter $s\in (0,1]$. This family is chosen in such a way
that $h(w;1)\equiv h(w)$ parametrizes the initial curve $\gamma $ whose
interior domain $D$ supports the eigenvalues. After the construction of $%
h(w;s)$ the support $D$ (resp. the harmonic moments $t_{j}$) also depends on 
$s$ under $s\mapsto D(s)$ (resp. $t\mapsto t_{j}(s)$). To state our result,
we denote by $\mathbf{\tau }=(\tau _{1},\ldots ,\tau _{n+1})$ a vector on
the affine space $Z\subset \mathbb{C}^{n+1}$ with $\tau _{1}=0$ and $%
\left\vert \tau _{2}\right\vert =1$.

\begin{Teorema}
\label{MainT} Consider the Elbau-Felder ensemble with $t_{0}>0$ and $%
s\mapsto \mathbf{t}(s)\in \mathbb{C}^{n+1}$ such that $t_{1}=0$, 
\begin{equation*}
t_{2}=\frac{\sqrt{1-s}}{2}\exp \left( is^{\Delta _{2}}\varphi \right) ~,
\end{equation*}%
and 
\begin{equation*}
t_{j}=s^{\Delta _{j}}\tau _{j}\qquad \text{for\qquad }3\leq j\leq n+1,
\end{equation*}%
with $\varphi \in \lbrack 0,2\pi )$, $\tau _{j}\in \mathbb{C}$, $\Delta
_{j}\geq 1$ and $s\in (0,1]$. There exist $r_{0}=r_{0}(\mathbf{\tau })>0$
such that for every $0<r<r_{0}$:

\begin{itemize}
\item[1)] There is a unique simple analytic closed polynomial curve $\gamma
=\gamma (s,r,\mathbf{\tau })$ of degree $n$, with external harmonic moments $%
\mathbf{t}(s)$ and area of interior domain $\pi t_{0}$ with $t_{0}$
depending on $\mathbf{\tau }$, $r$ and $s$.

\item[2)] The curve is parametrized by $h(w;s,r,\mathbf{\tau }%
)=rw+\sum_{j=0}^{n}r^{j}\alpha _{j}w^{-j}$, with $\left\vert w\right\vert =1$
and $\alpha _{j}=\alpha _{j}(s,\mathbf{\tau })$ is uniquely determined by $%
t_{0}$ and $\mathbf{t}(s)$.

\item[3)] The eigenvalue density is uniform within $D$, the interior domain
of $\gamma $, for every $s\in (0,1]$. Moreover, if $\Delta _{j}>1$ then the
Elbau-Felder ensemble can be conformally deformed, as $s$ goes to $0$, into
a Wigner ensemble with support on $[-2r,2r]$.
\end{itemize}
\end{Teorema}

\begin{Obs}
As long as $0<r<r_{0}$, $h(w;s,r,\mathbf{\tau })$ is a Riemann mapping from
the exterior of the unit disk onto the exterior domain $D$ of $\gamma $ and
the area $\pi t_{0}$ of the domain $D$ remains positive for all $s\in (0,1]$.
\end{Obs}

\begin{Obs}
The assumption $\left\vert t_{2}\right\vert <1/2$ in Theorem \ref%
{ElbauFelder} breaks down when the exterior domain $D_{-}(s)=\mathbb{C}%
\backslash D(s)$ is deformed into the slit domain $\mathbb{C}\backslash
\lbrack -2r,2r]$. In Section \ref{EHM} we generalize Elbau-Felder's results
using Crandall--Rabinowitz bifurcation theory from simple eigenvalues (see
e.g. \cite{Crandall}) to construct a parametrization that allow us to let $%
t_{2}\rightarrow 1/2$ maintaining the parameter $r$ away from $0$.
Elbau-Felder's parametrization, coming from the implicit function theorem
applied to the map $F$, mentioned right below Eq. (\ref{h}), defines a curve
in $\mathbb{R}\times \mathbb{C}^{n+1}$ that bifurcates at $t_{2}=1/2$.
\end{Obs}

\begin{Obs}
\label{3} If $\Delta _{j}>1$ then the $\alpha _{j}=\alpha _{j}(s,\mathbf{%
\tau })$ behave, for $s\rightarrow 0$, as

\begin{enumerate}
\item[(i)] $\alpha _{0}(s)=o(s)$,

\item[(ii)] $\alpha _{1}(s)=1-s/2$,

\item[(iii)] $\alpha _{j}(s)=o(s)$ for $1<j\leq n$.
\end{enumerate}

This relations will be used to prove item $\mathit{3)}$ of Theorem \ref%
{MainT}.
\end{Obs}

\begin{Obs}
Mashkov \textit{et al} (Sec. $6$ of \cite{EFensemble2}) considered a family
of simple closed analytic curves $\gamma (s)$ given by an equation $%
P(x,y/s)=0$, $w=x+iy$, converging to the segment of the real line $\left[
\alpha ,\beta \right] $ as $s\rightarrow 0$. At the level of formal
manipulations, they have shown that the eigenvalue density projected into $%
\left[ \alpha ,\beta \right] $ yields 
\begin{equation*}
\rho (x)=\lim_{s\rightarrow 0}\frac{\Delta y(x;s)}{s}=\sqrt{(x-\alpha
)(\beta -z)}M(x)
\end{equation*}%
where $\Delta y(x;s)$ is the width of $D(s)$ at the $x$ coordinate and $M(x)$
is a smooth function, regular at the edges.
\end{Obs}

This paper is organized as follows. Section \ref{Pre} presents some
preliminary results and introduces the two ingredients, the balayage problem
and the Schwarz function, required for the proof of Theorem \ref{MainT}. In
Section \ref{RM} we prove two auxiliary results, Propositions \ref{hn-I} and %
\ref{hn-II}. Section \ref{EHM} uses Crandall--Rabinowitz bifurcation theory
from simple eigenvalues to establish a smooth inverse map $F^{-1}$ in $%
\mathbb{R}\times \mathbb{C}^{n+1}$ about $t_{1}=0$ and $t_{2}=1/2$. Theorem %
\ref{BalV} in Section \ref{BM} gives an explicit expression of the Balayage
measure for the potential $V$. Section \ref{CD} concludes the proof of
Theorem \ref{MainT} based in Lemma \ref{MainLemma} and Section \ref{Ex}
gives some examples. We present in Section \ref{C} our conclusions and Lemma %
\ref{linear} is proved in Appendix \ref{A1}.

\section{Basic Setting}

\label{Pre}

\subsection{Eigenvalue Distribution for Normal Ensembles}

For normal unitarily invariant ensembles, we can write Eq. (\ref{Pm}) in
terms of the spectral coordinates. The joint probability of the eigenvalues $%
\{z_{i}\}_{i=1}^{N}\subset \Sigma $ of $M$ reads 
\begin{equation}
P_{N}(z_{1},\cdots ,z_{N})\propto \exp \Big\{{-\Big(2\sum_{1\leq i<j\leq
N}\log \left\vert z_{i}-z_{j}\right\vert ^{-1}+N\sum_{i=1}^{N}V(z_{i})\Big)%
\Big\}}.  \label{Pe3}
\end{equation}%
Introducing the empirical measure of the eigenvalues 
\begin{equation}
d\mu _{N}\left( z\right) =N^{-1}\sum_{i=1}^{N}\delta (z-z_{i})d^{2}z,
\label{medemp}
\end{equation}%
(\ref{Pe3}) can be written as 
\begin{equation*}
P_{N}(z_{1},\cdots ,z_{N})=Z_{N}^{-1}e^{-N^{2}I^{V}(\mu _{N})},
\end{equation*}%
where $Z_{N}$ is the normalization and 
\begin{equation}
I^{V}(\mu )\equiv \int \left( V(z)+U^{\mu }(z)\right) d\mu \left( z\right) 
\label{I}
\end{equation}%
is the total energy. The logarithmic potential associated with $\mu $ given
by 
\begin{equation}
U^{\mu }(z)\equiv \int \log \left\vert z-w\right\vert ^{-1}d\mu \left(
w\right) ~.  \label{U(z)}
\end{equation}%
The integrals with respect to (\ref{Pe3}) have, in the limit $N\rightarrow
\infty $, dominant contribution governed by a variational problem: 
\begin{equation}
E^{V}\equiv \underset{\mu \in \mathcal{M}(\Sigma )}{\inf }I^{V}(\mu )\text{
, \ \ }  \label{E}
\end{equation}%
where the infimum is taken over the set $\mathcal{M}(\Sigma )$ of Borel
probability measures in $\Sigma \subset \mathbb{C}$. If a probability
measure $\mu ^{V}$ satisfying 
\begin{equation*}
E^{V}=I(\mu ^{V})
\end{equation*}%
exists, it is called the equilibrium measure associated with $V$. The
empirical measure (\ref{medemp}) is known to converge weakly to a unique
equilibrium measure as $N\rightarrow \infty $ (see \cite{Deift} for
Hermitian ensembles and \cite{Hedenmalm} for normal ensembles).

\begin{Teorema}[Elbau-Felder]
\label{ElbauFelder} Consider the Elbau-Felder ensemble of degree $n$. There
is $\delta >0$ such that for all $0<t_{0}<\delta $ a unique equilibrium
measure $d\mu $ exists and is uniform within a domain $D\subset \Sigma $
that contains the origin: 
\begin{equation}
d\mu =\frac{1}{\pi t_{0}}\chi _{D}\left( z\right) d^{2}z;  \label{miEF}
\end{equation}%
$D$ is uniquely determined by the exterior harmonic moments ($t_{1}=0$) 
\begin{eqnarray}
\pi t_{0} &=&\int_{D}d^{2}z\qquad \,\mbox{the area of}\,D  \notag \\
t_{k} &=&\frac{-1}{\pi k}\int_{\mathbb{C}\backslash D}z^{-k}d^{2}z~,\text{
if }k=2,\cdots ,n+1,\text{ }  \notag \\
t_{k} &=&0\text{ if }k>n+1  \label{t0tk}
\end{eqnarray}%
and its boundary $\gamma $ is a simple closed analytic polynomial curve of
degree $n$; if $h(w)=rw+a_{0}+a_{1}/w+\cdots +a_{n}/w^{n}$, $|w|=1$,
parametrizes $\gamma $, then \textbf{\ } 
\begin{equation}
t_{0}=r^{2}-\sum_{j=1}^{n}j|a_{j}|^{2}.  \label{t0}
\end{equation}%
There exist homogeneous universal polynomials $P_{jk}\in \mathbb{Z}%
[r,a_{0},\cdots ,a_{k-j}]$ of degree $k-j+1$, $1\leq j\leq k\leq n+1$ such
that 
\begin{equation}
jt_{j}=\bar{a}_{j-1}r^{-j+1}+\sum_{k=j}^{n}\bar{a}_{k}r^{-k}P_{jk}(r,a_{0},%
\cdots ,a_{k-j})~  \label{tj-P}
\end{equation}%
is an invertible transformation from $\mathbb{R}\times \mathbb{C}^{n}$ into
itself in a neighborhood of $\left( r^{2},a_{0},a_{1}/r\ldots
,a_{n}/r^{n}\right) =(0,0,2\bar{t}_{2},\ldots ,(n+1)\bar{t}_{n+1})$. For
sufficiently small $r$, the function $h(w)$ is a Riemann mapping from the
exterior of the unit disk onto the exterior of $D$.
\end{Teorema}

The proof of existence of the equilibrium measure requires to verify
Euler--Lagrange type equations. Let%
\begin{equation}
E(z)=V(z)+\frac{2}{\pi t_{0}}\int_{D}\ln \left\vert \frac{z}{\zeta }%
-1\right\vert ^{-1}d^{2}\zeta   \label{Ez}
\end{equation}%
be a function defined in $\Sigma $ given by $V$ plus the logarithmic
potential (\ref{U(z)}) associated with the uniform measure $\mu $ in $D$.
Lemma 6.3 of \cite{ElbauFelder} shows that $E(z)=0$ holds for almost every $%
z\in D$. According to Corollary 3.5 of \cite{ElbauFelder}, $\mu $ is the
unique equilibrium measure if%
\begin{equation}
E(z)\geq 0\ \text{for every }z\in \Sigma \backslash D~.  \label{E+}
\end{equation}%
We extend Elbau--Felder's proof to the near--slit--domains in Section \ref%
{EHM}.

\subsection{Balayage Problem}

We tackle the problem of analyzing conformal deformations of normal ensemble
by means of balayage techniques \cite{Shapiro}. This allows us to solve the
problem focusing only on the behavior of the boundary $\gamma$ of $D$.
Therefore, in our approach the balayage technique plays a major role. Let $%
G\subset \overline{\mathbb{C}}$ be an open set and $\partial G$ its boundary.

Let $\nu $ be a probability measure on $G$ (such that $\nu \left( \overline{%
\mathbb{C}}\,\backslash G\right) =0$) and let the logarithmic potential $%
U^{\nu }$ (see Eq. (\ref{U(z)})) be finite and continuous on $G$. \textbf{%
The balayage problem} (or "sweeping out" problem) consists in finding a
probability measure $\widehat{\nu }$ with support on $\partial G$ such that 
\begin{equation}
U^{\nu }=U^{\widehat{\nu }}\text{ almost everywhere}\text{ on }\partial G.
\label{UU}
\end{equation}%
We call $\widehat{\nu }$ the balayage measure associated with $\nu .$
Throughout this work, we consider the following space of functions:

\begin{definition}
Let $G\subset \mathbb{C}$ be a bounded open set. We denote by $\mathcal{H}%
(G) $ the space of all holomorphic functions on $G$ and continuous on its
closure $\overline{G}$.
\end{definition}

If $G\subset \mathbb{C}$ is a bounded open set and $\nu $ is a probability
measure with compact support in $G$, then $\widehat{\nu }$ is the unique
measure supported in $\partial G$ satisfying (\ref{UU}) and such that $U^{%
\widehat{\nu }}(z)$ is bounded in $\partial G$. In addition, $\widehat{\nu }$
possesses the following property (see Theorem $II-4.1$ of \cite{Totik}): 
\begin{equation}
\int_{G}fd\nu =\int_{\partial G}fd\widehat{\nu }  \label{Bal}
\end{equation}%
holds for every $f\in \mathcal{H}(G)$.

We may choose $G$ the interior $\dot{D}=D\backslash \gamma $ of the compact
support $D$ of the equilibrium measure $\mu $ associated with $V$. In this
case, we write $\mu =\mu |_{\dot{D}}+\mu |_{\gamma }$ and sweep out only the
part $\mu |_{\dot{D}}$ lying on $G$: $\widehat{\mu }=\widehat{\mu |_{\dot{D}}%
}+\mu |_{\gamma }$. Since the equilibrium measure has no mass concentrated
in $\gamma $, we have $\mu |_{\gamma }\equiv 0$. The balayage measure
associated with the equilibrium measure $\mu $ is denoted simply by $%
\widehat{\mu }$ and has support in $\gamma $.

\subsection{Parametric Curves and Schwarz Function}

The following definitions will be important for the characterization of the
curves appearing in our main result. We shall start with the basic

\begin{definition}
A curve $\Gamma $ in $\mathbb{C}$ is said to be simple if there exist a
parametrization $t\mapsto h(t)$ for $t\in \lbrack a,b]$ such that $h(t)$ is
injective, i. e., if for all $x,y\in \lbrack a,b]$ $h(x)\not=h(y)$ when $%
x\not=y$. If $h(a)=h(b)$, in this case $\Gamma $ is said to be a simple
closed curve. A curve $\Gamma $ in $\mathbb{C}$ is said to be an analytic
curve if there exist a parametrization $t\mapsto h(t)$ for $t\in \lbrack a,b]
$ such that $h$ is analytic and $h^{\prime }(t)\not=0$ for $t\in \lbrack a,b]
$.
\end{definition}

Next we introduce the polynomial curves on the complex plane.

\begin{definition}
A curve $\Gamma $ in $\mathbb{C}$ is said to be a polynomial curve of degree
n if it is parametrically represented as 
\begin{equation}
h(w)=rw+a_{0}+a_{1}w^{-1}+\ldots +a_{n}w^{-n}.  \label{Rm}
\end{equation}%
with $r>0$, $a_{n}\not=0$ and $|w|=1$.
\end{definition}

We shall define the Schwarz function

\begin{definition}
Let $\Gamma $ in $\mathbb{C}$ be an analytic arc and let $\Omega $ be a
strip-like neighborhood of $\Gamma $. The Schwarz function $S$ is the unique
analytic function on $\Omega $ such that 
\begin{equation*}
S(z)=\bar{z}\ ,\qquad z\in \Gamma .
\end{equation*}
\end{definition}

For a treatise on the Schwarz function with applications see \cite%
{Davis,Shapiro}.

\begin{Obs}
Hereafter, $\gamma$ denotes a simple closed analytic polynomial curve.
Moreover, $S$ stands for the Schwarz function of $\gamma$.
\end{Obs}

Schwarz function $S$ will play a major role in the conformal deformation of
the Elbau-Felder ensemble. We shall show the balayage measure is
proportional to the Schwarz function $S$.

\section{Riemann Map\label{RM}}

We shall prove some auxiliary results, which concern the behavior of the
family $h(w;s)$ as the $\gamma (s)$ is deformed to the real line.

\begin{Propo}
\label{hn-I} Let $h(w;s)=rw+\displaystyle\sum_{j=0}^{n}a_{j}(s)w^{-j}$, $%
\left\vert w\right\vert =1$, be for each $s\in (0,1]$ the parametrization of
a closed polynomial curve $\gamma (s)$ of degree $n$ with $h(\cdot
;1)=\gamma $ a simple curve. Let the condition%
\begin{equation}
\xi (s):=r-\sum_{j=1}^{n}j\left\vert a_{j}(s)\right\vert >0  \label{tau}
\end{equation}%
be satisfied for every $s\in (0,1]$. Then, for each $s\in (0,1]$, the $%
\gamma (s)$ remains a simple polynomial curves and $h(w;s)$, as a map from
the exterior of the unit disk into the exterior of $\gamma (s)$, is
biholomorphic (a Riemann map). Furthermore, for every $s\in \lbrack \delta
,1]$, $0<\delta <1$, $t_{0}(s)$ is bounded away from zero.
\end{Propo}

\noindent \textit{Proof.} Let us begin with the estimation of $t_{0}$. It
follows from (\ref{tau}) that $r>|a_{j}(s)|$ holds for every $j$.
Multiplying $\xi (s)$ by $r$ it yields that 
\begin{equation*}
0<r^{2}-\sum_{j=1}^{n}j\left\vert a_{j}(s)\right\vert
r<r^{2}-\sum_{j=1}^{n}j\left\vert a_{j}(s)\right\vert ^{2}=t_{0}(s)
\end{equation*}%
is bounded away from zero for $s\in \lbrack \delta ,1]$, $0<\delta <1$. Eq. (%
\ref{tau}) also implies that $h(w;s)$ is an analytic curve, that is, the
derivative of $h(w;s)$ with respect to $w$, denoted by $h^{\prime }(w;s)$ is
bounded away from zero:%
\begin{equation*}
\left\vert h^{\prime }(w;s)\right\vert =\left\vert
r-\sum_{j=1}^{n}ja_{j}(s)w^{-(j+1)}\right\vert \geq
r-\sum_{j=1}^{n}j|a_{j}(s)||w|^{-(j+1)},
\end{equation*}%
where in the last passage we have used the triangular inequality. Since we
are analyzing the exterior of the unity circle $\left\vert w\right\vert >1$,
we have 
\begin{equation}
\left\vert h^{\prime }(w;s)\right\vert \geq r-\sum_{j=1}^{n}j\left\vert
a_{j}(s)\right\vert >0\text{ }\forall s\in (0,1],  \label{hw}
\end{equation}%
which also holds in a small neighborhood of $\left\vert w\right\vert =1$.

Now, for every $w$, $z\in \mathbb{C}$ with $\left\vert w\right\vert
=\left\vert z\right\vert =1$, by the triangular inequality, 
\begin{eqnarray}
\left\vert h(w;s)-h(z;s)\right\vert &\geq &r\left\vert w-z\right\vert
-\sum_{j=1}^{n}\left\vert a_{j}(s)\right\vert \left\vert \frac{1}{w^{j}}-%
\frac{1}{z^{j}}\right\vert  \notag \\
&=&r\left\vert w-z\right\vert -\sum_{j=1}^{n}\left\vert a_{j}(s)\right\vert
\left\vert w^{j}-z^{j}\right\vert  \notag \\
&\geq &\xi (s)\left\vert w-z\right\vert >0  \label{wz}
\end{eqnarray}%
if $w\neq z$. The last passage follows from $\left\vert
w^{j}-z^{j}\right\vert \leq j\left\vert w-z\right\vert $ which can be shown
using the telescopic identity%
\begin{equation*}
w^{j}-z^{j}=w^{j-1}(w-z)+w^{j-2}(w-z)z+\cdots +(w-z)z^{j-1}
\end{equation*}%
together with the triangular inequality. Equations (\ref{wz}) and (\ref{hw})
imply that the map $h(\cdot ;s):S^{1}\longrightarrow \mathbb{C}$ is an
embedding, $\gamma (s)$ is a simple curve and $h(w;s)$ is a Riemann map from
the exterior of the unit circle onto the exterior of $\gamma (s)$ for every $%
s\in (0,1]$. The polynomial curve $\gamma (s)$ with $0<s<1$ preserves all
properties assumed for $\gamma (1)=\gamma $, concluding the proof
Proposition \ref{hn-I}. $\Box $

\begin{Obs}
\label{appeal}Proposition \ref{hn-I} has an intuitive appeal. To a
polynomial curve $\gamma (s)$ fail to be simple it has to develop a cusp.
However, when $h(w_{c},s)$ form a cusp, we have $h^{\prime }(w_{c},s)=0$,
this situation is prevented as long as $\xi (s) > 0$. Proposition \ref{hn-I}
gives a sufficient condition for $h^{\prime }(w,s)\neq 0$ and show that (\ref%
{tau}) is also sufficient for $\gamma (s)$ to remain a simple curve.
\end{Obs}

The next result concerns the conditions $(i-iii)$ of \textit{Remark} \ref{3}
and the deformation $\gamma (s)$ to the real line as $s\rightarrow 0$.

\begin{Propo}
\label{hn-II} Consider a polynomial curve $\gamma (s)$ of degree $n$
parametrized by $h(w;s)$ satisfying conditions $(i-iii)$. Then, $%
\lim_{s\rightarrow 0}h(w;s)$ maps the exterior of the unit disk onto $%
\mathbb{C}\backslash \lbrack -2r,2r]$.
\end{Propo}

\noindent \textit{Proof.} As $h(w,s)$ is a Riemann map for all $s\in (0,1]$
by Proposition \ref{hn-I}, then it suffices to show that $\lim_{s\rightarrow
0}\gamma (s)=[-2r,2r]$. Indeed, since $\left\vert w\right\vert =1$ we may
choose $w=e^{i\theta }$ with $\theta \in \lbrack 0,2\pi ]$. The Riemann map
reads 
\begin{equation*}
h(e^{i\theta };s)=re^{i\theta }+a_{0}(s)+a_{1}(s)e^{-i\theta
}+\sum_{j=2}^{n}a_{j}(s)e^{-ij\theta },
\end{equation*}%
which, under the hypotheses, yields 
\begin{equation*}
h(w;s)=r(e^{i\theta }+e^{-i\theta })+o(s),
\end{equation*}%
implying that $\lim_{s\rightarrow 0}h(e^{i\theta },s)=2r\cos \theta \in
\lbrack -2r,2r]$ for all $\theta \in \lbrack 0,2\pi ]$. $\Box $

\section{Exterior Harmonic Moments of Near-to-Slit Domains\label{EHM}}

Let $\mathbf{t}=(t_{1},t_{2},\ldots ,t_{n+1})$ be the exterior harmonic
moment of the domain $D$ -- containing the origin and bounded by $\gamma $
-- and let $\pi t_{0}$ be the area of $D$.

When a given collection $\mathbf{t}$ of $n+1$ complex numbers, together with 
$t_{0}>0$, determines a simple polynomial curve $\gamma $ of degree $n$? We
refer to Theorem 5.3 of \cite{ElbauFelder} for the solution to this moment
problem. If $t_{1}=0$ and $\left\vert t_{2}\right\vert <1/2$, then every
complex numbers $t_{2},\ldots ,t_{n+1}$ determine a curve $\gamma $ with
these properties provided $t_{0}$ is sufficiently small.  We shall give a
proof of this result in a language more appropriated for the generalization
needed.

The map $\left( \rho ,\mathbf{\alpha }\right) \in \mathbb{R}_{+}\times 
\mathbb{C}^{n+1}\longmapsto \left( t_{0},\mathbf{t}\right) \in \mathbb{R}%
_{+}\times \mathbb{C}^{n+1}$ defined by (\ref{t0}) and by the contour
integral (\ref{tj}),%
\begin{eqnarray}
jt_{j} &=&\frac{1}{2\pi i}\oint_{|w|=1}\bar{h}(w^{-1})h^{\prime
}(w)h^{-j}(w)dw~  \notag \\
&=&\sum_{k=j-1}^{n}\bar{\alpha}_{k}\text{Res}\left( w^{k-j}\frac{%
1-\sum_{l=1}^{n}l\alpha _{l}\rho ^{l}/w^{l+1}}{\left( 1+\sum_{l=0}^{n}\alpha
_{l}\rho ^{l}/w^{l+1}\right) ^{j}};\infty \right)  \label{tj-h}
\end{eqnarray}%
taking residues at infinity: 
\begin{eqnarray}
t_{0} &=&\rho -\sum\nolimits_{j=1}^{n}j|\alpha _{j}|^{2}\rho ^{j}  \notag \\
t_{j} &=&\frac{1}{j}\bar{\alpha}_{j-1}-\bar{\alpha}_{j}\alpha _{0}-\left( 1+%
\frac{1}{j}\right) \alpha _{1}\bar{\alpha}_{j+1}\rho +O\left( \rho
^{2}\right) ,\ \ \ \ 1\leq j\leq n+1  \label{ttt}
\end{eqnarray}%
with $\rho =r^{2}$, $\alpha _{j}=r^{-j}a_{j}$, $0\leq j\leq n$, and $\alpha
_{j}=0$ if $j>n$, has a smooth inverse in a neighborhood of $%
(0,0,t_{2},\ldots ,t_{n+1})$ provided $\left\vert t_{2}\right\vert \neq 1/2$.

In this section, we show how the inverse function theorem is applied in this
situation and extend it for the case $\left\vert t_{2}\right\vert =1/2$. We
also verify whether the inverse determines a simple polynomial curve $\gamma
=\partial D$ and whether a measure $\mu $, uniform in a near-to-slit domain $%
D$, is the equilibrium measure of the Elbau-Felder ensemble.

\subsection{Bifurcating Curves}

We observe that $\rho ^{\ast }=\alpha _{0}^{\ast }=0$ and $\alpha _{j}^{\ast
}=(j+1)\bar{t}_{j+1}$, $j=1,\ldots ,n$ solve the equations (\ref{ttt}) for $%
\left( \rho ,\alpha _{0},\ldots ,\alpha _{n}\right) $ and because it takes
complex conjugation of the $t_{j+1}$, we look (\ref{ttt}) as a map from $%
\mathbb{R}_{+}\times \mathbb{C}^{n+1}\times \mathbb{C}^{n+1}$ into itself 
\begin{equation}
F:(\rho ,\mathbf{\alpha },\mathbf{\bar{\alpha}})\longmapsto \left( t_{0},%
\mathbf{t},\mathbf{\bar{t}}\right)  \label{F}
\end{equation}%
and we write%
\begin{eqnarray}
0 &=&\frac{F(\rho ,\mathbf{\alpha }^{\ast }+\rho \mathbf{\varphi },\mathbf{%
\bar{\alpha}}^{\ast }+\rho \mathbf{\bar{\varphi}})-\left( t_{0},\mathbf{t},%
\mathbf{\bar{t}}\right) }{\rho }  \notag \\
&=&l(1,\mathbf{\varphi },\mathbf{\bar{\varphi}})+p(\rho ,\mathbf{\varphi },%
\mathbf{\bar{\varphi}})-\left( \tau _{0},\mathbf{0},\mathbf{0}\right) ~.
\label{FLr}
\end{eqnarray}%
Here, $F$ has been expanded in Taylor series about $(\rho ^{\ast },\mathbf{%
\alpha }^{\ast },\mathbf{\bar{\alpha}}^{\ast })$ with remainder $\rho p(\rho
,\mathbf{\varphi },\mathbf{\bar{\varphi}})=O\left( \rho ^{2}\right) $, $\tau
_{0}=t_{0}/\rho =O\left( 1\right) $ and $l$ is the linear map 
\begin{equation}
l(1,\mathbf{\varphi },\mathbf{\bar{\varphi}})=L\left[ \rho ^{\ast },\mathbf{%
\alpha }^{\ast },\mathbf{\bar{\alpha}}^{\ast }\right] \left( 
\begin{array}{c}
1 \\ 
\mathbf{\varphi } \\ 
\mathbf{\bar{\varphi}}%
\end{array}%
\right) =\left( 
\begin{array}{ccc}
1-4\left\vert t_{2}\right\vert ^{2} & \mathbf{0}^{T} & \mathbf{0}^{T} \\ 
-\mathbf{\bar{v}} & -\bar{K} & J^{-1} \\ 
-\mathbf{v} & J^{-1} & -K%
\end{array}%
\right) \left( 
\begin{array}{c}
1 \\ 
\mathbf{\varphi } \\ 
\mathbf{\bar{\varphi}}%
\end{array}%
\right)  \label{l}
\end{equation}%
where $\mathbf{0}$ is the null column vector in $\mathbb{C}^{n+1}$, $\bar{v}%
_{j}=2\left( 1+1/j\right) (j+2)\bar{t}_{2}t_{j+2}$ if $1\leq j<n$ with $%
v_{n}=v_{n+1}=0$, $O$ denotes the $(n+1)\times (n+1)$ null matrix, $J=\text{%
diag}\left\{ j\right\} _{j=1}^{n+1}$ and $K$ is a $(n+1)\times (n+1) $
matrix with $\bar{K}_{i1}=(i+1)t_{i+1}$ for $1\leq i\leq n$ and $0$
otherwise.

Since $L$ is invertible for $\left\vert t_{2}\right\vert \neq 1/2$, (\ref{F}%
) has a smooth inverse defined in a neighborhood of $(t_{0},\mathbf{t}%
)=(0,0,t_{2},\ldots ,t_{n+1})$. The implicit function theorem applied to (%
\ref{FLr}) (with $\mathbf{t}\in \mathbb{C}^{n+1}$ fixed) uniquely defines
two smooth curves parametrized by $\rho $: 
\begin{eqnarray*}
\mathbf{\varphi }(\rho ) &=&T\mathbf{v}+B\mathbf{\bar{v}}+O(\rho ) \\
\tau _{0}(\rho ) &=&1-4\left\vert t_{2}\right\vert ^{2}+O(\rho )
\end{eqnarray*}%
on $\mathbb{C}^{n+1}$ and $\mathbb{R}_{+}$, respectively, where $%
B=(1-4\left\vert t_{2}\right\vert ^{2})^{-1}JK$ and $T=J+2\bar{t}_{2}B$. We
note that the leading order in $\rho $ of equation (\ref{FLr}) can be
written as%
\begin{eqnarray}
1-4\left\vert t_{2}\right\vert ^{2} &=&\tau _{0}~  \notag \\
\left( 
\begin{array}{cc}
-\bar{K} & J^{-1} \\ 
J^{-1} & -K%
\end{array}%
\right) \left( 
\begin{array}{c}
\mathbf{\varphi } \\ 
\mathbf{\bar{\varphi}}%
\end{array}%
\right) &=&\left( 
\begin{array}{c}
\mathbf{\bar{v}} \\ 
\mathbf{v}%
\end{array}%
\right)  \label{sist1}
\end{eqnarray}%
whose solution is given in Appendix \ref{A1}.

The function $t_{0}=\rho \tau _{0}(\rho )$ is monotone increasing in $\rho
\in \left[ 0,a\right] $ for $\left\vert t_{2}\right\vert <1/2$ and
sufficiently small $a$, and its inverse, $\rho (t_{0},\mathbf{t})$, is a
well defined function of $t_{0}$ and $\mathbf{t}$. The inverse of (\ref{ttt}%
) in $\mathbb{R}_{+}\times \mathbb{C}^{n+1}$ thus reads%
\begin{equation}
(t_{0},\mathbf{t})\longmapsto \left( \rho (t_{0},\mathbf{t}),\mathbf{\alpha }%
^{\ast }+\rho (t_{0},\mathbf{t})~\mathbf{\varphi \circ }\rho (t_{0},\mathbf{t%
})\right) ~.  \label{inv}
\end{equation}

The above application of the implicit function theorem breaks down if the
second harmonic moment $t_{2}$ tends to $1/2$ and this is the case when the
external domain of $\gamma (s)$ tends, as $s\rightarrow 0$, to the slit
domain $\mathbb{C}/\left[ -2r,2r\right] $ (see Proposition \ref{hn-II}). We
need, therefore, to extend Theorem 5.3 of \cite{ElbauFelder} to include this
case. For this, we shall rescale all components of $\mathbf{t}$, excepted $%
t_{2}$ whose modulus square will be denoted by $\lambda =\left\vert
t_{2}\right\vert ^{2}$, as%
\begin{equation}
t_{j}=(1-4\lambda )\tau _{j}~,\ j\neq 2  \label{ttal}
\end{equation}%
and apply the constructive bifurcation theory from a simple eigenvalue
developed by Crandall and Rabinowitz (see e.g. \cite{Crandall}). The theory
applied to equation (\ref{FLr}) uses the implicit function theorem with the
role of $\rho $ and $\lambda $ exchanged. Instead of the two parametric
curves $\varphi _{j}=\varphi _{j}(\rho )$ and $\tau _{0}=\tau _{0}(\rho )$,
we consider $\tilde{\varphi}_{j}=\tilde{\varphi}_{j}(\rho ,\lambda )$ and $%
\tilde{\tau}_{0}=\tilde{\tau}_{0}(\rho ,\lambda )$ as a function of $\lambda 
$ for $\left( \rho ,\mathbf{\tau }\right) \in \mathbb{C}^{n+1}$ fixed, where 
$\mathbf{\tau }=\left( \tau _{1},\tau _{2},\ldots ,\tau _{n+1}\right) $ is a
vector in the affine space of $\mathbb{C}^{n+1}$ with $\left\vert \tau
_{2}\right\vert =1$, denoted by $Z$.

We observe that $\left( \rho ^{\ast },\mathbf{\alpha }^{\ast }\right) $ with 
$\left\vert t_{2}\right\vert =1/2$ is a bifurcation point for (\ref{ttt})
since every neighborhood of this point contains a solution which differs
from (\ref{inv}). Note that the tangent map of (\ref{ttt}) at $\left( \rho
^{\ast },\mathbf{\alpha }^{\ast }\right) $ with $\left\vert t_{2}\right\vert
=1/2$ is singular i. e., $L$ is not invertible at $\lambda =\left\vert
t_{2}\right\vert ^{2}=1/4$. Using bifurcation theory, we shall construct a
pair of smooth curves for fixed $\left( \rho ,\mathbf{\tau }\right) \in 
\mathbb{C}^{n+1}$: $\tilde{\varphi}_{j}=\tilde{\varphi}_{j}(\lambda )$ and $%
\tilde{\tau}_{0}=\tilde{\tau}_{0}(\lambda )$, $\lambda \in (1/4-b,1/4]$ for
some $b>0$, such that $\tilde{\varphi}_{j}(1/4)=\varphi _{j}(\rho
=0,\left\vert t_{2}\right\vert =1/2)$ and $\tilde{\tau}_{0}(1/4)=\tau
_{0}(\rho =0,\left\vert t_{2}\right\vert =1/2)$.

\begin{Propo}
\label{inverse}Given $\left( \rho ,\mathbf{\tau }\right) \in \mathbb{R}%
_{+}\times Z\simeq \mathbb{C}^{n+1}$, there exist two uniquely defined
smooth curves, $\lambda \mapsto \mathbf{\zeta }(\rho ,\lambda ,\mathbf{\tau }%
)$ on $\mathbb{C}^{n+1}$ and $\lambda \mapsto \tilde{\tau}_{0}(\rho ,\lambda
,\mathbf{\tau })$ on $\mathbb{R}$, defined by (\ref{zeta}), such that the
inverse of (\ref{ttt}) in a neighborhood of $(0,0,\mathbf{\tau })$ in $%
\mathbb{R}_{+}\times \mathbb{R}_{+}\times Z$ is written as%
\begin{equation}
(t_{0},\rho ,\mathbf{\tau })\longmapsto \left( \rho ,\mathbf{\alpha }^{\ast
}\circ \lambda (t_{0},\rho ,\mathbf{\tau })+\rho (\mathbf{\alpha }_{0}+%
\mathbf{\zeta }\circ \lambda (t_{0},\rho ,\mathbf{\tau }))\right) ~
\label{inv1}
\end{equation}%
where $\lambda (t_{0},\rho ,\mathbf{\tau })$ is a function of $t_{0},\rho $
and $\mathbf{\tau }$ which is the unique solution for $\lambda $ of $%
t_{0}=\rho \tilde{\tau}_{0}(\rho ,\lambda ,\mathbf{\tau })$ in the domain $%
0\leq \lambda \leq 1/4$, $\rho \in \left[ 0,\hat{r}^{2}\right] $ with $\hat{r%
}=\hat{r}(\mathbf{\tau })$ sufficiently small. Moreover, there exist $\bar{r}%
=\bar{r}(\mathbf{\tau })>0$ such that $h(w)$, defined by (\ref{h}) with $r=%
\sqrt{\rho }$ and $a_{j}=r\alpha _{j}$ given by (\ref{inv1}), parametrizes a
simple polynomial curve of order $n$ which can be deformed to the segment $%
\left[ -2r,2r\right] $, for every fixed $r<\bar{r}$.
\end{Propo}

\noindent \textit{Proof.} It suffices, for the first statement, to verify
the hypothesis of Theorem 1 in Sec. 3 of \cite{Crandall}. For this, it is
convenient to write (\ref{FLr}) as 
\begin{equation}
l_{0}(1,\mathbf{\varphi },\mathbf{\bar{\varphi}})+(1-4\lambda )l_{1}(1,%
\mathbf{\varphi },\mathbf{\bar{\varphi}})+\tilde{p}(\lambda ,\rho ,\mathbf{%
\varphi },\mathbf{\bar{\varphi}})-\left( \tau _{0},\mathbf{0},\mathbf{0}%
\right)   \label{FLLr}
\end{equation}%
where $l_{0}$ and $l_{1}$ are linear maps in $\mathbb{R}_{+}\times \mathbb{C}%
^{n+1}\times \mathbb{C}^{n+1}$, with $l_{0}=$ $\left. l\right\vert _{\lambda
=1/4}$ and $4l_{1}=-\left. \partial (l+p)/\partial \lambda \right\vert
_{\lambda =1/4}$:%
\begin{equation}
l_{0}(1,\mathbf{\varphi },\mathbf{\bar{\varphi}})=L_{0}\left( 
\begin{array}{c}
1 \\ 
\mathbf{\varphi } \\ 
\mathbf{\bar{\varphi}}%
\end{array}%
\right) =\left( 
\begin{array}{ccc}
0 & \mathbf{0}^{T} & \mathbf{0}^{T} \\ 
\mathbf{0} & -\bar{K}_{0} & J^{-1} \\ 
\mathbf{0} & J^{-1} & -K_{0}%
\end{array}%
\right) \left( 
\begin{array}{c}
1 \\ 
\mathbf{\varphi } \\ 
\mathbf{\bar{\varphi}}%
\end{array}%
\right)   \label{l0}
\end{equation}%
\begin{equation}
l_{1}(1,\mathbf{\varphi },\mathbf{\bar{\varphi}})=L_{1}\left( 
\begin{array}{c}
1 \\ 
\mathbf{\varphi } \\ 
\mathbf{\bar{\varphi}}%
\end{array}%
\right) =\left( 
\begin{array}{ccc}
1 & \mathbf{0}^{T} & \mathbf{0}^{T} \\ 
-\mathbf{\bar{v}}_{1} & -\bar{K}_{1} & O \\ 
-\mathbf{v}_{1} & O & -K_{1}%
\end{array}%
\right) \left( 
\begin{array}{c}
1 \\ 
\mathbf{\varphi } \\ 
\mathbf{\bar{\varphi}}%
\end{array}%
\right)   \label{l1}
\end{equation}%
with $\left( \bar{\mathbf{v}}_{1}\right) _{j}=\left( 1+1/j\right) (j+2)\tau
_{j+2}/\tau _{2}+O(\rho )$, $\left( K_{0}\right) _{ij}=\tau _{2}\delta
_{i1}\delta _{j1}+O(\rho )$ and, given 
\begin{equation*}
\bar{\mathbf{w}}_{1}=\left( -\tau _{2}/2,3\tau _{3},\ldots ,(n+1)\tau
_{n+1},0\right) \,,
\end{equation*}%
$K_{1}$ is, up to the leading order in $\rho $, a $(n+1)\times (n+1)$ matrix
with $\mathbf{w}_{1}$ in its first column and $0$ everywhere else. Moreover, 
$\tilde{p}(\lambda ,\rho ,\mathbf{\varphi },\mathbf{\bar{\varphi}})$ is a
smooth map from $\mathbb{R}_{+}\times \mathbb{R}_{+}\times \mathbb{C}%
^{n+1}\times \mathbb{C}^{n+1}$ to $\mathbb{R}_{+}\times \mathbb{C}%
^{n+1}\times \mathbb{C}^{n+1}$ which satisfies%
\begin{equation*}
\tilde{p}(1/4,\rho ,\mathbf{\varphi },\mathbf{\bar{\varphi}})=\tilde{p}%
(\lambda ,0,\mathbf{\varphi },\mathbf{\bar{\varphi}})=\frac{\partial \tilde{p%
}}{\partial \lambda }(1/4,\rho ,\mathbf{\varphi },\mathbf{\bar{\varphi}})=0~.
\end{equation*}

We observe that (\ref{FLr}), and consequently (\ref{FLLr}), extends by
continuity to $\rho =0$. In this way, we define%
\begin{equation*}
G(\rho ,\lambda ,\mathbf{\zeta ,\bar{\zeta}})=\frac{F(\rho ,\mathbf{\alpha }%
^{\ast }+\rho (\mathbf{\alpha }_{0}+\mathbf{\zeta }),\mathbf{\bar{\alpha}}%
^{\ast }+\rho (\mathbf{\bar{\alpha}}_{0}+\mathbf{\bar{\zeta}}))-\left( t_{0},%
\mathbf{t},\mathbf{\bar{t}}\right) }{\rho }~,
\end{equation*}%
for $(\lambda ,\rho ,\mathbf{\zeta ,\bar{\zeta}})\in \left[ 1/4-b,1/4+b%
\right] \times \left[ 0,a\right] \times \mathbb{C}^{n+1}\times \mathbb{C}%
^{n+1}$ where, for $\mathbf{\alpha }_{0}=(\sqrt{\tau _{2}},0,\ldots ,0)$, $%
(1,\mathbf{\alpha }_{0},\mathbf{\bar{\alpha}}_{0})$ is an eigenvector of $%
L_{0}$ associated with the (simple) null eigenvalue. It follows by (\ref%
{FLLr})-(\ref{l1}), that 
\begin{equation*}
G(\rho ,\lambda ,\mathbf{\zeta ,\bar{\zeta}})=l_{0}(0,\mathbf{\zeta },%
\mathbf{\bar{\zeta}})+(1-4\lambda )l_{1}(1,\mathbf{\alpha }_{0}+\mathbf{%
\zeta },\mathbf{\bar{\alpha}}_{0}+\mathbf{\bar{\zeta}})+\tilde{p}(\lambda
,\rho ,\mathbf{\alpha }_{0}+\mathbf{\zeta },\mathbf{\bar{\alpha}}_{0}+%
\mathbf{\bar{\zeta}})-\left( \tau _{0},\mathbf{0},\mathbf{0}\right) ~.
\end{equation*}%
and the implicit function theory can be applied to $G=0$ provided the
derivative of $G(\rho ,\lambda ,\mathbf{\zeta ,\bar{\zeta}})$ about $%
(\lambda ,\mathbf{\zeta ,\bar{\zeta}})=\left( 1/4,\mathbf{0},\mathbf{0}%
\right) $, defined by the linear map%
\begin{equation*}
(\lambda ,\mathbf{\zeta ,\bar{\zeta}})\longmapsto l_{0}(0,\mathbf{\zeta },%
\mathbf{\bar{\zeta}})-4\lambda l_{1}(1,\mathbf{\alpha }_{0},\mathbf{\bar{%
\alpha}}_{0})=\left( 
\begin{array}{ccc}
-4 & \mathbf{0}^{T} & \mathbf{0}^{T} \\ 
4(\mathbf{\bar{v}}_{1}+\mathbf{\bar{w}}_{1}\sqrt{\tau _{2}}) & -\bar{K}_{0}
& J^{-1} \\ 
4(\mathbf{v}_{1}+\mathbf{w}_{1}/\sqrt{\tau _{2}}) & J^{-1} & -K_{0}%
\end{array}%
\right) \left( 
\begin{array}{c}
\lambda \\ 
\mathbf{\zeta } \\ 
\mathbf{\bar{\zeta}}%
\end{array}%
\right) ~,
\end{equation*}%
is nonsingular (see Theorem 1 in Sec. 3 of \cite{Crandall}). Since it is
always invertible, $G=0$ (with $\rho >0$ and $\mathbf{\tau }\in Z$ fixed)
uniquely defines two smooth parametric curves: 
\begin{eqnarray}
\mathbf{\zeta }(\lambda ) &=&\bar{\tau}_{2}K_{0}(\mathbf{v}_{1}+\mathbf{w}%
_{1}/\sqrt{\tau _{2}})+K_{0}(\mathbf{\bar{v}}_{1}+\mathbf{\bar{w}}_{1}\sqrt{%
\tau _{2}})+O(1-4\lambda ,\rho )  \notag \\
\tilde{\tau}_{0}(\lambda ) &=&1-4\lambda +O(\left( 1-4\lambda \right)
^{2},\rho )  \label{zeta}
\end{eqnarray}%
on $\mathbb{C}^{n+1}$ and on $\mathbb{R}_{+}$. As in the previous case,
writing $\eta =1-4\lambda $, the leading order of equation $G=0$ reads%
\begin{eqnarray}
\eta &=&\tilde{\tau}_{0}  \notag \\
\left( 
\begin{array}{cc}
-\bar{K}_{0}-\eta \bar{K}_{1} & J^{-1} \\ 
J^{-1} & -K_{0}-\eta K_{1}%
\end{array}%
\right) \left( 
\begin{array}{l}
\mathbf{\zeta } \\ 
\mathbf{\bar{\zeta}}%
\end{array}%
\right) &=&\eta \left( 
\begin{array}{c}
\mathbf{\bar{v}}_{1}+\mathbf{\bar{w}}_{1}\sqrt{\tau _{2}} \\ 
\mathbf{v}_{1}+\mathbf{w}_{1}/\sqrt{\tau _{2}}%
\end{array}%
\right)  \label{sist2}
\end{eqnarray}%
from where (\ref{zeta}) can be readily obtained.

The function $t_{0}=\rho \tilde{\tau}_{0}(\rho ,\lambda )$ is monotone
decreasing in $\lambda $, for $0<\lambda \leq 1/4$, $\rho \in \left[ 0,\hat{r%
}^{2}\right] $ with $\hat{r}$ sufficiently small, and its inverse, $\lambda
(t_{0},\rho ,\mathbf{\tau })$, is a well defined function of $t_{0},\rho $
and $\mathbf{\tau }$. The inverse of (\ref{ttt}) in $\mathbb{R}_{+}\times 
\mathbb{R}_{+}\times \hat{Z}$ is thus given by (\ref{inv1}). By (\ref{alphaj}%
), $\hat{r}=\hat{r}(\mathbf{\tau })$ is defined by 
\begin{equation*}
2\sum_{j=1}^{n}j(j+1)^{2}\left\vert \tau _{j+1}\right\vert ^{2}\hat{r}%
^{2(j-1)}=1~.
\end{equation*}

To conclude the proof, it remains to verify that 
\begin{equation}
h(w)=rw+\alpha _{0}+r\frac{\alpha _{1}}{w}+\cdots +r^{n}\frac{\alpha _{n}}{%
w^{n}}\ ,~\left\vert w\right\vert =1  \label{hwr}
\end{equation}%
with 
\begin{equation}
\alpha _{j}=\alpha _{j}^{\ast }\circ \lambda (t_{0},\rho ,\mathbf{\tau }%
)+\rho (\alpha _{0,j}+\zeta _{j}\circ \lambda (t_{0},\rho ,\mathbf{\tau }))
\label{alphaj}
\end{equation}%
parametrizes a simple polynomial curve $\gamma =\gamma (t_{0})$. For this,
in order $\gamma $ to approach the segment $I=\left[ -2r,2r\right] $, as the
area of its interior $\pi t_{0}$ tends to $0$, we set $\tau _{2}=1$ and let $%
\tau _{3}$ and $\rho =r^{2}$ so that 
\begin{equation}
\alpha _{0}=\rho \left( \alpha _{0,1}+\zeta _{1}\right) =12\Re \mathrm{e}%
(\tau _{3})r^{2}+O(r^{4},t_{0})  \label{alpha0}
\end{equation}%
remains inside $I$, closed to the origin. Now, it is enough to verify the
hypothesis (\ref{tau}) of Proposition \ref{hn-I}. It follows immediately
from 
\begin{eqnarray}
\alpha _{1} &=&\left( 1-\frac{t_{0}}{2r^{2}}\right) +O(t_{0}^{2})  \notag \\
\alpha _{j} &=&\frac{t_{0}}{r^{2}}(j+1)\bar{\tau}_{j+1}+O(t_{0}^{2})~,\ \
2\leq j\leq n~  \label{alpha}
\end{eqnarray}%
for every $\mathbf{\tau }$, that 
\begin{equation*}
\xi =r-\sum_{j=1}^{n}jr^{j}\left\vert \alpha _{j}\right\vert \geq \frac{t_{0}%
}{r^{2}}\left( \frac{r}{2}-\sum_{j=2}^{n}j(j+1)r^{j}\left\vert \tau
_{j+1}\right\vert \right) >0
\end{equation*}%
holds for $r<\bar{r}$ where $\bar{r}$ is defined by equating the expression
between parenthesis to $0$. $\Box $ 

\begin{Obs}
$r_{0}$ in Theorem \ref{MainT} is thus given by $\min \left( \bar{r},\hat{r}%
\right) $. 
\end{Obs}

\begin{Obs}
The results of Proposition \ref{inverse} can be adapted for (\ref{ttal})
scaled with different power of $\left( 1-4\lambda \right) $: $%
t_{j}=(1-4\lambda )^{\Delta _{j}}\tau _{j}$, $\Delta _{j}\geq 1$ for each $%
j\neq 2$. 
\end{Obs}

Note that Proposition \ref{inverse} holds for simple closed analytic curves
provided $l $, $l_{0}$ and $l_{1}$ are Frech\'{e}t derivatives of $F$ with
respect to an appropriate Banach space.

\subsection{Equilibrium Measure}

Equation (\ref{E+}) is equivalent to 
\begin{equation}
0\leq \mathcal{E}(w):=E(h(w))=\frac{1}{t_{0}}\left( \left\vert
h(w)\right\vert ^{2}-\left\vert h(1)\right\vert ^{2}+2\Re \mathrm{e}%
\int_{1}^{w}\bar{h}(\zeta ^{-1})h^{\prime }(\zeta )d\zeta \right)  \label{EE}
\end{equation}%
for $w\in h^{-1}\left( \Sigma \right) $ such that $\left\vert w\right\vert
\geq 1$ (see eq. (16) of \cite{ElbauFelder}), and it suffices to verify only
for $\left\vert w\right\vert \geq 1/R$, where $R=\max \left\{ \left\vert
w\right\vert :h^{\prime }(w)=0~,\ w\in \mathbb{C}\right\} $ is the critical
radius of $\gamma $, and $r$ sufficiently small.

Let $h^{(0)}(w)=rw+\alpha _{0}+r\alpha _{1}/w$ be the parametrization of an
ellipse that approximates $\gamma $ and let $(t_{0}^{(0)},t_{2}^{(0)})$
denote the corresponding external harmonic moment. We denote by $\mathcal{E}%
^{(0)}(w)$ the function defined in (\ref{EE}) for the ellipse and note that,
by Subsection 6.2 of \cite{ElbauFelder}, $t_{0}^{(0)}\mathcal{E}^{(0)}(w)$
remains bounded for $w\in h^{-1}\left( \Sigma \right) $ such that $%
\left\vert w\right\vert \geq 1$. For $(1-t_{0}/r^{2})^{-1}\leq \left\vert
w\right\vert <r^{-\alpha }$, $0<\alpha <1/3$, we have 
\begin{eqnarray*}
\frac{1}{t_{0}}\left( h(w^{-1})-h^{(0)}(w^{-1})\right) h^{\prime }(w) &=&%
\frac{1}{t_{0}}\frac{r^{2}\alpha _{2}}{w^{-2}}r+O(t_{0})=O(r^{1-2\alpha }) \\
\frac{1}{t_{0}}h^{(0)}(w^{-1})\left( h^{\prime }(w)-h^{(0)\prime }(w)\right)
&=&\frac{1}{t_{0}}r\alpha _{1}w\frac{-2r^{2}\alpha _{2}}{w^{3}}+O(t_{0})=O(r)
\\
\frac{t_{0}-t_{0}^{(0)}}{t_{0}^{2}} &=&\frac{1}{t_{0}^{2}}r^{4}\alpha
_{2}^{2}+O(t_{0})=O(1) \\
t_{0}\mathcal{E}^{(0)}(w) &=&O(r^{2-\alpha })
\end{eqnarray*}%
uniformly in $t_{0}$ as $t_{0}\rightarrow 0$. Consequently, 
\begin{eqnarray*}
\left\vert \mathcal{E}(w)\right\vert &\geq &\frac{1}{t_{0}}\left\vert t_{0}%
\mathcal{E}^{(0)}(w)\right\vert -\left\vert \mathcal{E}(w)-\mathcal{E}%
^{(0)}(w)\right\vert \\
&\geq &\frac{1}{t_{0}}\left\vert t_{0}\mathcal{E}^{(0)}(w)\right\vert -\frac{%
1}{t_{0}}\left\vert \left\vert h(w)\right\vert ^{2}-\left\vert
h^{(0)}(w)\right\vert ^{2}-\left\vert h(1)\right\vert ^{2}+\left\vert
h^{(0)}(1)\right\vert ^{2}\right\vert \\
&&-\frac{2}{t_{0}}\left\vert \Re \mathrm{e}\int_{1}^{w}\left( h(\zeta
^{-1})-h^{(0)}(\zeta ^{-1})\right) h^{\prime }(\zeta )d\zeta \right\vert \\
&&-\frac{2}{t_{0}}\left\vert \Re \mathrm{e}\int_{1}^{w}h(\zeta ^{-1})\left(
h^{\prime }(\zeta )-h^{(0)\prime }(\zeta )\right) d\zeta \right\vert
-\left\vert \frac{t_{0}-t_{0}^{(0)}}{t_{0}^{2}}t_{0}\mathcal{E}%
^{(0)}(w)\right\vert
\end{eqnarray*}%
is strictly positive for $t_{0}$ sufficiently small. For $w\in h^{-1}\left(
\Sigma \right) $ such that $\left\vert w\right\vert \geq r^{-\alpha }$ we
may proceed exactly as in \cite{ElbauFelder}.

\section{Explicit Balayage Measures\label{BM}}

To explicitly obtain the balayage measure in terms of the potentials we need
an auxiliary result.

\begin{Teorema}
\label{dmiu} If $V:\Sigma \rightarrow \mathbb{R}$ is a potential defined in
a compact set $\Sigma \subset \mathbb{C}$ with continuous second partial
derivatives in its interior, then the variational problem (\ref{E}) is
attained at a unique equilibrium measure $\mu ^{V}$ supported in a compact
set $D\subset \Sigma $ given by 
\begin{equation}
d\mu ^{V}(z)=\frac{1}{4\pi }\Delta Vd^{2}z,  \label{distmi}
\end{equation}
at almost every (w.r.t. the Lebesgue measure $d^{2}z$) interior point of $D$
, where $\Delta =\partial ^{2}/\partial x^{2}+\partial ^{2}/\partial y^{2}$
is the Laplace operator.
\end{Teorema}

\noindent \textit{Proof.} Using the smoothness of $V$ the proof follows from
Theorems $I-1.3$ pp. 27 and $II-1.3$ pp. 85 of \cite{Totik}. $\Box $

We derive an explicit equation for the balayage measure in terms of the
potential $V$. The result will be of fundamental importance to establish the
conformal deformation of the Elbau-Felder ensemble.

A direct application of this theorem to the potential $V$ shows that the
density of eigenvalues is indeed uniform within $D$. The Elbau-Felder
potential reads $V(z)=(z\bar{z}-p(z)-\bar{p(z)})/t_{0}$, remembering that $%
\Delta =4\partial _{z}\partial _{\bar{z}}$ and $p(z)$ is analytic ($%
\partial_{\bar{z}} p(z)=\partial _{z}\overline{p(z)}=0$), it follows that $%
d\mu (z)=1/(\pi t_{0})d^{2}z$ within $D$. Moreover, we have the following
result

\begin{Teorema}
\label{BalV} Let $V$ be an Elbau--Felder potential. The balayage measure $%
\widehat{\mu }$ associated with the equilibrium measure $\mu ^{V}$ with
support on $\gamma =\partial D$ is 
\begin{equation}
d\widehat{\mu }(z)=\frac{1}{2\pi it_{0}}S(z)dz,  \label{bal}
\end{equation}
where $dz$ is the measure of arc length on $\gamma $.
\end{Teorema}

\noindent \textit{Proof.} The main ingredient is the Green theorem. Our
hypotheses on $V$, ${D}$ and $\partial {D}$ guarantee that the Green theorem
is applicable. Using that $\Delta =4\partial _{z\bar{z}}$ and $f\in \mathcal{%
H}(D)$, by (\ref{distmi}) and applying the Green theorem we have 
\begin{eqnarray}
\int_{D}f(z)d\mu ^{V}(z) &=&\frac{1}{\pi }\int_{D}\partial _{\overline{z}%
}\left( f(z)\partial _{z}V(z)\right) d^{2}z  \notag \\
& = & \frac{1}{2i\pi }\int_{\partial D}f(z) \partial_{z}V(z) dz  \notag
\end{eqnarray}
and by the balayage measure property Eq. (\ref{Bal}), we have 
\begin{equation*}
\int_{\partial D}f(z)d\widehat{\mu }^{V}(z)=\frac{1}{2i\pi }\int_{\partial
D}f(z)\partial _{z}V(z)dz
\end{equation*}
holds for every continuous function on $\partial D$, from where it follows
the continuity of the balayage measure with respect to the Lebesgue measure.
Next, we have that 
\begin{equation*}
d\hat{\mu}=\frac{1}{2\pi i}\partial _{z}V(z)=\frac{1}{t_{0}}\left( \overline{%
z}- \sum_{k=2}^{n+1}kt_{k}z^{k-1}\right)
\end{equation*}
the term of the sum does not contribute to a contour integral with respect
to $\widehat{\mu }$ for test functions $f\in \mathcal{H}(D)$ (by continuity
and by the Cauchy theorem). Then, by using the definition $\bar{z}=S(z)$, $%
z\in \gamma $, of the Schwarz function we obtain Eq. (\ref{bal}), concluding
the proof. $\Box $

Eq. (\ref{bal}) relates the equilibrium measure with the Schwarz function of
boundary of the support. We shall conclude, by applying the Cauchy theorem,
that only the branch cut in the interior domain $D$ of the Schwarz function
will contribute to the balayage measure. Thus, questions concerning the
equilibrium measure $\mu $ turns to the behavior of the Schwarz function. It
turns out that, except when $\gamma $ is a line or a circle arc, the Schwarz
function $S$ always has a branch cut \cite{Davis}. Our next result draws
some conclusions on the behavior of the branch cuts.

\begin{Propo}
If $\gamma $ is a simple closed analytic curve, then the Schwarz function $S$
associated with $\gamma $ must have branching points in its interior. \label%
{1}
\end{Propo}

\noindent \textit{Proof.} Take $f\in \mathcal{H}(D)$ such that $%
\int_{D}f(z)d\mu \not=0$. By the property (\ref{Bal}) of balayage measure
together with the Theorem \ref{BalV}, we have 
\begin{equation*}
\int_{D}f(z)d\mu (z)=\frac{1}{2\pi it_{0}}\oint_{\gamma }f(z)S(z)dz.
\end{equation*}
Suppose now that $S$ has no branch point in the interior of $\gamma $, then
by using the Cauchy theorem we would then conclude, contrarily to the
hypothesis, that $\int_{D}f(z)d\mu =0$. Therefore, we must have an even
number of branch points inside $D$. Another important result is that the
branch point never touches the curve $\gamma $. This can be proved by noting
that $\gamma $ is, by hypothesis, a simple analytic polynomial curve and the
Schwarz function $S$ must be analytic on $\gamma $ and on its neighborhood
(see \cite{ElbauFelder,Davis}), what exclude the case of the branch point
touching the curve $\gamma $, concluding the proof. $\Box $

\section{Conformal Deformation\label{CD}}

In this section we conclude the proof of our main result, Theorem \ref{MainT}%
. Since the Riemann map $h(z;s)$ is conformal from the exterior of the unit
disk onto the exterior of $\gamma (s)$ it has a well defined inverse from
the exterior of $\gamma (s)$ onto the exterior of the unit disk. We shall
denote its inverse by $H(z;s)$: 
\begin{equation}
h(H(z;s);s)=z\,  \label{hInv}
\end{equation}
for all $z\,$in the exterior of $\gamma (s)$ and $s\in (0,1]$.

The Schwarz functions $S$ can be related to the Riemann map $h$ and its
inverse $H$ by the following:

\begin{Propo}
\label{Sch}Let $\gamma $\ be a polynomial curve parametrized by $h$. Then
the Schwarz function is a biholomorphic map in a neighborhood of $\gamma $\
and is given by 
\begin{equation*}
S(z)=\bar{h}\left( \frac{1}{H(z)}\right) ,
\end{equation*}
\ where $\bar{h}(w)=rw+\bar{a}_{0}+\bar{a}_{1}w^{-1}+\cdots +\bar{a}
_{n}w^{-n}$.
\end{Propo}

The proof can be found in Refs. \cite{ElbauFelder,Davis}.

\begin{Propo}
\label{Hn} If $h(z;s)$ satisfies the conditions $(i-iii)$ of Remark \ref{3},
then the inverse function $H(z;s)$ of $h(z;s)$ reads 
\begin{equation}
H(z;s)=\frac{z+\sqrt{z^{2}-4ra_{1}(s)}}{2r}+o(s)  \label{H}
\end{equation}%
for $0<s<\varepsilon $ and $\varepsilon $ small enough.
\end{Propo}

\noindent \textit{Proof.} We have by Eqs. (\ref{Rm}) and (\ref{hInv}) 
\begin{equation*}
h(H(z);s)=rH(z;s)+\sum_{j=0}^{n}a_{j}(s)H(z;s)^{-j}=z
\end{equation*}%
which yields 
\begin{equation}
rH^{2}(z;s)-(z+o(s))H(z;s)+a_{1}(s)=0~.  \label{rootH}
\end{equation}

Solving the equation for $H$ the result follows. Note that $\left\vert
H\right\vert =\mathcal{O}(1)$ in the neighborhood of $\gamma (s)$. $\Box $

\begin{Obs}
Equation (\ref{H}) has branches, namely the minus square root, which can be
directly seen from Eq. (\ref{rootH}). We do not consider the minus square
root, because it does not yield $S(z)=\bar{z}$ on $\gamma(s) $.
\end{Obs}

\begin{Lema}
\label{MainLemma} Let $d\hat{\mu}(z;s)$ be the balayage measure supported on 
$\gamma (s)$ of the Elbau-Felder ensemble. Assume the conditions \textbf{(i}%
--\textbf{iii)} to hold. Then, given $\varepsilon $ sufficiently small, for $%
0<s<\varepsilon $ and $f\in \mathcal{H}(D(s))$, we have 
\begin{equation}
\oint_{\gamma (s)}f(z)d\hat{\mu}(z;s)=\frac{1}{2\pi r^{2}}\int_{-2r}^{2r}f(z)%
\sqrt{4r^{2}-x^{2}}dx+o\left( 1\right) ~.  \label{cml}
\end{equation}
\end{Lema}

\noindent \textit{Proof. } Given $z$ in a neighborhood of $\gamma (s)$ and $%
\varepsilon >0$ sufficiently small, the Schwarz function of $\gamma (s)$ for 
$0<s<\varepsilon $ by Proposition \ref{Sch} reads 
\begin{equation*}
S(z;s)=r\frac{1}{H(z;s)}+\bar{a}_{1}(s)H(z;s)+\ldots + \bar{a}_{n}(s)H^n(z;s)
\end{equation*}
which together with (\ref{H}), yields 
\begin{equation}
S(z;s)=Ez+\Lambda \sqrt{z^{2}-4ra_{1}(s)}+g(z;s)  \label{Sn(z;e)}
\end{equation}
where $E= (r^{2}+ | a_{1}^{2}(s) | )/ 2ra_{1}(s)$, $\Lambda = (\left\vert
a_{1}^{2}(s)\right\vert -r^{2})/2ra_{1}(s)$, and $g(z;s)$ is $o(s)$. By (\ref%
{bal}), we have 
\begin{equation*}
I(s)=\oint_{\gamma (s)}f(z)d\hat{\mu}=\oint_{\gamma (s)}f(z)\frac{S(z;s)}{
2\pi it_{0}(s)}dz,
\end{equation*}
for $f\in \mathcal{H}(D(s))$, which together with (\ref{Sn(z;e)}), gives 
\begin{equation}
I(s)=\oint_{\gamma (s)}\frac{f(z)}{2\pi it_{0}(s)}\Big(Ez+\Lambda \sqrt{
z^{2}-4ra_{1}(s)}+g(z;s)\Big)dz.  \label{Meio}
\end{equation}
\noindent The linear term in (\ref{Meio}) will give no contribution by
Cauchy theorem. Now we need to estimate term depending on $g$. Since $g$ is
analytic in a neighborhood of $\gamma(s)$, we may have the bound 
\begin{equation}
\frac{1}{2\pi t_0(s)} \left\vert \oint_{\gamma (s)}g(z;s)dz\right\vert \leq 
\frac{ \mbox{ max}_{z \in \gamma(s)} |g(z;s)|}{2 \pi t_0(s)} \oint_{\gamma
(s)}dz.  \label{intg}
\end{equation}
Since arc length $\oint_{\gamma (s)}dz$ is finite, $\mbox{ max}_{z \in
\gamma(s)} |g(z;s)| = o(s)$, and by (\ref{t0}) and conditions $(ii)$ of 
\textit{Remark } \ref{3}. 
\begin{equation*}
t_{0}(s)=r^{2}-|a_{1}(s)|^{2}+o(s^{2})= s+o(s),
\end{equation*}
we have that the r.h.s. of (\ref{intg}) is $o(s)/t_0(s) = o(1)$. Hence, $%
\Lambda/ t_0(s) = 1/2r^2 + \mathcal{O}(s)$. Also, for $z\in \gamma (s)$, 
\begin{equation*}
\frac{\sqrt{4ra_{1}(s)-z^{2}}}{ra_{1}(s)}=\frac{\sqrt{4r^{2}-z^{2}}}{r^{2}} +%
\mathcal{O}(s).
\end{equation*}
Combining the estimates we have Therefore, 
\begin{equation}
I(s)=\frac{-1}{4\pi }\oint_{\gamma (s)}f(z)\frac{\sqrt{4r^2-z^2}}{r^2}dz+o(1)
\label{Ieq}
\end{equation}
By using Proposition (\ref{hn-II}), $\gamma(s)$ may be deformed until it
coincides with the branch cut of the square root. Because the branch the
integrand along the around the cut will not be the same. In the first term
running from $2r$ to $-2r$, the square root becomes $\sqrt{4r^{2}-x^{2}}$,
while in the second running from $-2r$ to $2r$ the square root becomes $-%
\sqrt{4r^{2}-x^{2}}$. This leads to 
\begin{equation*}
I(s)=\frac{1}{2\pi r^{2}}\int_{-2r}^{2r}f(x)\sqrt{4r^{2}-x^{2}}dx + o(1)
\end{equation*}
concluding the proof. $\Box $

\subsection{Proof of Theorem \protect\ref{MainT}:}

We have constructed a one parameter family of curves $\gamma (s)$
parametrized by $h(w;s)$ with $\left\vert w\right\vert =1$. We have shown,
in Proposition \ref{hn-I}, that under the hypothesis (\ref{tau}) the family $%
\gamma (s)$ is composed of simple closed analytic polynomial curves of
degree $n$. Moreover, $h(w;s)$ acts as the Riemann map from the exterior of
the unity circle onto the exterior of $\gamma (s)$ for $s\in (0,1]$. By
Proposition \ref{inverse}, the problem of determining the exterior moments $%
t_{j}$ out of simple closed analytic polynomial curves has a unique solution
also when $t_{2}$ tends to $1/2$.

The balayage techniques enables us to sweep all the eigenvalues to the
boundary $\gamma (s)$ of the support $D(s)$ and analyze the deformation only
focusing on $\gamma (s)$. In Theorem \ref{BalV} we deduce an explicit
equation relating the balayage measure with the potential $V$ for $f \in 
\mathcal{H}$. It turns out that the balayage measure is proportional to the
Schwarz function. Our focus changes to the behavior of the Schwarz functions
for polynomial curves.

In Proposition \ref{1} we show that the Schwarz function of every simple
closed polynomial curve must have a branch cut which lies in the interior of
the curve and the branch points never touch the curve itself. The limit $%
s\rightarrow 0$, the curve $\gamma $ is smashed into the line, see
Proposition \ref{hn-I}. The same limit is analyzed in Lemma \ref{1} where we
show that the Schwarz function converges to a function of the type $\sqrt{
r^{2}-x^{2}}$.

This implies that as $s\rightarrow 0$ the normal ensemble converges to a
Hermitian one (all the eigenvalues are real). Using these results and
comparing the left hand side of (\ref{cml}) with the semicircle law (Eq. \ref%
{miW}), we conclude that conformal deformations of Elbau-Felder ensembles to
the real line yields a Wigner ensemble, as claimed. $\Box $

\section{Examples \label{Ex}}

\subsection{ Potential of degree 2}

Taking $n=1$ in (\ref{pzn}), the potential $V$ reads 
\begin{equation*}
V\left( z\right) =\frac{z\overline{z}-t_{2}z^{2}-\overline{t}_{2}\overline{z}%
^{2}}{t_{0}}.
\end{equation*}

Consider the one parameter family for Riemann maps $h(w;s)=rw+a_{1}(s)w^{-1}$
in which, without loss of generality, $r$ and $a_{1}(s)=r(1-s)$ for $s\in
(0,1]$. Note that this parametrization fulfills the conditions $(i-iii)$ by
taking $s\mapsto s/2$. This yields $h(w;s)=r(w+(1-s)w^{-1})~$ $%
t_{0}=r^{2}(2s-s^{2})$ and $2t_{2}=1-s$, by relations (\ref{t0}) and (\ref%
{tj-P}). The support of the eigenvalues is then given by 
\begin{equation*}
D(s)=\left\{ x+iy\in \mathbb{C}:\frac{x^{2}}{r^{2}(2-s)^{2}}+\frac{y^{2}}{%
r^{2}s^{2}}\leq 1\right\} .
\end{equation*}%
with major and minor semi--axis given, respectively, by $(2-s)$ and $s$. For
every $s\in (0,1]$ the Elbau-Felder conditions is satisfied, namely, $%
t_{1}(s)=0$ and $\left\vert t_{2}(s)\right\vert <1/2$. Condition (\ref{tau}%
): $\xi (s)=r-a_{2}(s)=s>0$ for $s\in (0,1]$, holds and conditions $(i--iii)$
are trivially satisfied.

The equilibrium measure associated to this potential is uniform within the
support $D(s)$ and may be continuous deformed up to the limit $%
\lim_{s\rightarrow 0}h(w;s)=2\Re e(w)$ for $\left\vert w\right\vert =1$, in
which its boundary becomes the line segment $\left[ -2r,2r\right] $.

To show that the balayage measure converges to the Wigner measure, we
consider the Schwarz function, which can be computed using the formulae in
Ref. \cite{Davis} or using $H(w;s)$ together with Proposition \ref{Sch}. By
equation (\ref{bal}) and the Cauchy theorem, only the non-holomorphic part
of the Schwarz function gives a contribution to the integrals. For test
functions $f\in \mathcal{H} (D(s))$ and the properties of the balayage
measure (\ref{Bal}), we have 
\begin{eqnarray*}
\int_{D(s)}f(z)d\mu (z;s) &=& \frac{1}{2\pi r^2 \left( 1-s\right) }\int_{-2
r \sqrt{1-s}}^{2 r \sqrt{1-s}}f(x)\sqrt{4r^2 (1 - s) -x^{2}}dx
\end{eqnarray*}

It is also important to mention that during the whole process of conformal
deformation of the normal ensemble to the Hermitian, the logarithm potential
along the curve is kept unchanged, being an invariant in the process of
deformation. On can see this by computing the integral explicitly.

\subsection{Potentials of Degree $3$: Breaking Down the Hypotheses}

If we break the conditions $(i-iii)$ of $Remark$ \ref{3} ($\Delta_j > 1$),
then it is impossible to deform the normal ensemble into a Hermitian one
keeping, at same time, the regularity of the curve $\gamma (s)$ and a non
trivial support. We shall illustrate this scenario for $n=2$ taking $%
t_1=t_2=0$. The potential then reads 
\begin{equation*}
V\left( z\right) =\frac{z\overline{z}-t_{3}z^{3}-\overline{t}_{3}\overline{z}%
^{3}}{t_{0}}.
\end{equation*}
By (\ref{t0}) we have $t_{0}(s)=r^{2}-2|a_{2}(s)|^{2}$ and $%
3t_{3}(s)=a_{2}(s)/r^{2}$. The Riemann map is written as $%
h(w;s)=rw+a_{2}(s)w^{-2}$. Note that $a_{1}=0$ violating condition $(ii)$.
Therefore, the deformation of $\gamma (s)$ to a segment of the real line is
impossible. The regularity of $\gamma (s)$, which is guaranteed by condition
(\ref{tau}), requires $|a_{2}|<r/2$, while the positive area condition on $%
t_{0}$ requires $|a_{2}|<r/\sqrt{2}$.

There is no parametrization that keeps the regularity of $\gamma $. When the
area $\pi t_{0}$ of $D(s)$ converges to zero ($a_{2}\rightarrow r/2$), $%
\gamma (s)$ develops cusps and is no longer regular. Keeping the regularity
of the curve, the limit $t_{0}\rightarrow 0$ is possible only taking $%
r\rightarrow 0$. Therefore, as the area converges to zero $\gamma $
collapses to a point.

\section{Conclusion \label{C}}

We study the conformal deformations of Elbau-Felder ensembles to Hermitian
ensembles. Special attention is paid to Wigner ensembles, for which the
density of eigenvalues follows the semicircular law. The result presented in
Theorem \ref{MainT} is in a way universal, in the sense that it does not
depend on the initial ensemble one starts with. Two ingredients are used to
prove the main result. We use bifurcation theory from a single eigenvalue
for extending the moment problem to near slit domains and the concept of
balayage measure and Schwarz function for assisting the convergence of
normal ensembles to Wigner ensembles.

It would be interesting to analyze to deformation of the Elbau-Felder
ensembles to the real line without the assumption on the convergence rate.
In this case the resulting measure, if well defined, might deviate from the
semicircular law, and it might also depend on the initial ensemble under
consideration.

We would like to acknowledge the financial of the Brazilian agencies FAPESP
under grant 07/04579-2 (TP) and CNPq under 151976/2007-4 (AMV). DHUM is
partially supported by CNPq. \appendix

\section{Solution of (\protect\ref{sist1})\label{A1}}

\begin{Lema}
\label{linear}The linear system of equations (\ref{sist1}) has a unique
solution $\mathbf{\varphi }=T\mathbf{v}+B\mathbf{\bar{v}}$ where $%
B=(1-\left\vert k\right\vert ^{2})^{-1}JK$ with $k=K_{11}$ and $T=J+\bar{k}B$%
, provided $\left\vert k\right\vert \neq 1$.
\end{Lema}

\noindent \textit{Proof.} Equation (\ref{sist1}) is solvable if and only if 
\begin{equation*}
\left( 
\begin{array}{cc}
B & T \\ 
\bar{T} & \bar{B}%
\end{array}%
\right) \left( 
\begin{array}{cc}
-\bar{K} & J^{-1} \\ 
J^{-1} & -K%
\end{array}%
\right) =\left( 
\begin{array}{cc}
I & 0 \\ 
0 & I%
\end{array}%
\right)
\end{equation*}%
holds for some $(n+1)\times (n+1)$ complex matrices $B$ and $T$ and this is
equivalent to%
\begin{eqnarray}
-B\bar{K}+TJ^{-1} &=&I  \notag \\
BJ^{-1}-TK &=&0~.  \label{BK}
\end{eqnarray}%
Let us assume that $B$ has its first column given by $\mathbf{b}$ and $0$
everywhere else (so $B$ has the same form of $K$) and let $T=J+\bar{k}B$.
Note that $KA=a_{11}K$ ($BA=a_{11}B$) holds for any matrix $A=[a_{ij}]$.
Substituting $T$ in (\ref{BK}), we have%
\begin{eqnarray*}
BJ^{-1}-TK &=&BJ^{-1}-(J+\bar{k}B)K \\
&=&B-JK-\left\vert k\right\vert ^{2}B=0
\end{eqnarray*}%
which implies%
\begin{eqnarray*}
B &=&\frac{1}{1-\left\vert k\right\vert ^{2}}JK \\
T &=&J+\frac{\bar{k}}{1-\left\vert k\right\vert ^{2}}JK~.
\end{eqnarray*}%
and concludes the proof of the lemma. The uniqueness follows by linearity.$%
~~~\Box $

\end{document}